\documentclass{amsart}
\usepackage{amsmath, amssymb}
\usepackage{amsfonts}
\usepackage[arrow,matrix,curve,cmtip,ps]{xy}

\usepackage{amsthm}

\allowdisplaybreaks

\newtheorem{theorem}{Theorem}[section]
\newtheorem{lemma}[theorem]{Lemma}

\newtheorem{corollary}[theorem]{Corollary}

\newtheorem*{theorem*}{Theorem}
\theoremstyle{remark}
\newtheorem{remark}[theorem]{Remark}
\newtheorem{definition}[theorem]{Definition}
\newtheorem{example}[theorem]{Example}

\numberwithin{equation}{section}

\newcommand{\N}{\mathbb{N}}
\newcommand{\R}{\mathbb{R}}
\newcommand{\C}{\mathbb{C}}

\newcommand{\K}{\mathcal{K}}

\begin{document}
\title[Stability of $C^*$-algebras associated to graphs]{Stability of
$\boldsymbol{C^*}$-algebras associated to graphs}

\author{Mark Tomforde 
}

\address{Department of Mathematics\\ Dartmouth College\\
Hanover\\ NH 03755-3551\\ USA}

\curraddr{Department of Mathematics\\ University of Iowa\\
Iowa City\\ IA 52242\\ USA}

\email{tomforde@math.uiowa.edu}


\date{\today}

\begin{abstract}
We characterize stability of graph $C^*$-algebras by giving five
conditions equivalent to their stability.  We also show that if $G$ is a
graph with no sources, then $C^*(G)$ is stable if and only if each vertex in
$G$ can be reached by an infinite number of vertices.  We use this
characterization to realize the stabilization of a graph $C^*$-algebra. 
Specifically, if $G$ is a graph and $\tilde{G}$ is the graph formed by adding
a head to each vertex of $G$, then $C^*(\tilde{G})$ is the stabilization of
$C^*(G)$; that is, $C^*(\tilde{G}) \cong C^*(G) \otimes \K$.
\end{abstract}

\maketitle

\section{Introduction}

In 1980 Cuntz and Krieger introduced a class of $C^*$-algebras generated by
families of partial isometries satisfying relations determined by a finite
matrix with entries in $\{0,1\}$.  These Cuntz-Krieger algebras were
initially studied because of their appearance in the study of topological
Markov chains.  Later it was found that they also have important
parallels with certain kinds of dynamical systems (e.g.~shifts of finite
type).

Since their inception, Cuntz-Krieger algebras have been generalized in an
extraordinary number of ways.  One generalization whose study has proven
particularly fruitful are the $C^*$-algebras associated to directed graphs. 
In 1982 Watatani noted that one could view the Cuntz-Krieger algebra
associated to a finite matrix $A$ as the $C^*$-algebra associated to the
finite directed graph with adjacency matrix $A$ \cite{Wat}.  However,
these ideas were not more fully explored until the late 1990's when Kumjian,
Pask, Raeburn, and Renault \cite{KPRR} introduced $C^*$-algebras associated to
locally finite graphs (i.e.~possibly infinite graphs in which each vertex
emits and receives a finite number of edges).  Not long after this it was
shown in \cite{BPRS} that many of the same results also hold for
$C^*$-algebras associated to row-finite graphs (i.e.~possibly infinite graphs
in which every vertex emits finitely many edges) and often the same
techniques can be applied to prove these results.  In the early 2000's
$C^*$-algebras associated to arbitrary directed graphs were finally
considered \cite{FLR}.  Unlike the generalization from locally finite to
row-finite graphs, it was found that extending results to $C^*$-algebras of
arbitrary graphs often involved significant modifications to statements of
theorems as well as the development of new techniques for their proofs.

In this paper we consider the notion of stability for $C^*$-algebras
associated to arbitrary directed graphs.  Recall that a $C^*$-algebra $A$
is said to be stable if $A \cong A \otimes \K$, where $\K$ denotes the
compact operators on a separable infinite-dimensional Hilbert space. 
Furthermore, if $A$ is a $C^*$-algebra then one may form its
stabilization $A \otimes \K$.  Since $\K \otimes \K \cong \K$, one has
that the stabilization of a $C^*$-algebra is stable.

If $G$ is a graph and $C^*(G)$ is its associated $C^*$-algebra, then we
prove in Theorem~\ref{stable-thm} that the stability of $C^*(G)$ is
equivalent to five other conditions.  This theorem generalizes a result of
Hjelmborg \cite[Theorem~2.14]{Hje2} in which a characterization of stability
for $C^*$-algebras of locally finite graphs was obtained.  Both of these
results make use of a nontrivial characterization of stability due to R\o
rdam and Hjelmborg \cite{HR}.  However, our proof of Theorem~\ref{stable-thm}
will involve techniques significantly different from Hjelmborg's proof of
\cite[Theorem~2.14]{Hje2}.  Furthermore, in addition to applying to
$C^*$-algebras of arbitrary graphs, Theorem~\ref{stable-thm} is different from
\cite[Theorem~2.14]{Hje2} in another respect, namely that it includes a
characterization in terms of the graph traces on all of $G$, rather than on a
special subgraph of $G$ as in Condition~(d) of \cite[Theorem~2.14]{Hje2}.

In Corollary~\ref{left-inf-stable} we show that there is a particularly nice
characterization of stability for $C^*(G)$ when $G$ has no sources: 
If $G$ is a graph with no sources, then $C^*(G)$ is stable if and only if
every vertex of $G$ can be reached by an infinite number of vertices.  This
gives an easily verifiable condition for determining the stability of the
$C^*$-algebra solely in terms of the graph.

Building off this characterization, in \S\ref{stabilization-sec} we develop a
method for realizing the stabilization of a graph algebra and we show that
it is also a graph algebra.  If $G$ is a graph, then we obtain
a new graph $\tilde{G}$ by adding a ``head"
$$
\xymatrix{
\cdots \ar[r] & \bullet \ar[r] & \bullet \ar[r] & \bullet \ar[r] & v\\ }
$$
to each vertex $v$ in $G$.  We prove in Theorem~\ref{stabilization-gr-alg-thm}
that $C^*(\tilde{G})$ is the stabilization of $C^*(G)$; that is,
$C^*(\tilde{G}) \cong C^*(G) \otimes \K$.  As a corollary we have that the
class of graph algebras is closed under stabilization.

\section{Preliminaries}

We provide some basic facts about graph algebras and refer
the reader to \cite{KPR}, \cite{BPRS}, and \cite{BHRS} for more
details.  A (directed) graph $G=(G^0, G^1, r, s)$ consists of a
countable set $G^0$ of vertices, a countable set $G^1$ of edges,
and maps $r,s:G^1 \rightarrow G^0$ identifying the range and
source of each edge.  A vertex $v \in G^0$ is called a
\emph{sink} if $|s^{-1}(v)|=0$, and $v$ is called an
\emph{infinite emitter} if $|s^{-1}(v)|=\infty$.  If $v$ is either a sink or
an infinite emitter, then we call $v$ a \emph{singular vertex}.  A graph
$G$ is said to be \emph{row-finite} if it has no infinite emitters.

If $G$ is a graph we define a \emph{Cuntz-Krieger $G$-family}
to be a set of mutually orthogonal projections $\{p_v : v \in
G^0\}$ and a set of partial isometries $\{s_e : e \in G^1\}$ with
orthogonal ranges which satisfy the \emph{Cuntz-Krieger relations}:
\begin{enumerate}
\item $s_e^* s_e = p_{r(e)}$ for every $e \in G^1$;
\item $s_e s_e^* \leq p_{s(e)}$ for every $e \in G^1$;
\item $p_v = \sum_{s(e)=v} s_e s_e^*$ for every $v
\in G^0$ that is not a singular vertex.
\end{enumerate}
The \emph{graph algebra $C^*(G)$} is defined to be the
$C^*$-algebra generated by a universal Cuntz-Krieger $G$-family.  

A \emph{path} in $G$ is a sequence of edges $\alpha =
\alpha_1 \alpha_2 \ldots \alpha_n$ with $r(\alpha_i) =
s(\alpha_{i+1})$ for $1 \leq i < n$, and we say that $\alpha$ has
length $|\alpha| = n$.  We let $G^n$ denote the set of all paths
of length $n$, and we let $G^* := \bigcup_{n=0}^\infty G^n$ denote
the set of finite paths in $G$.  Note that vertices are
considered paths of length zero.  The maps $r,s$ extend to
$G^*$, and for $v,w \in G^0$ we write $v \geq w$ if there exists a path
$\alpha \in G^*$ with $s(\alpha)=v$ and $r(\alpha) = w$.  Also for a path
$\alpha := \alpha_1 \ldots \alpha_n$ we define $s_\alpha := s_{\alpha_1}
\ldots s_{\alpha_n}$.  It is a consequence of the Cuntz-Krieger relations
that $C^*(G) = \overline{\textrm{span}} \{ s_\alpha s_\beta^* : \alpha, \beta
\in G^* \text{ and } r(\alpha) = r(\beta)\}$.

We say that a path $\alpha := \alpha_1 \ldots \alpha_n$ of length $1$ or
greater is a \emph{loop} if $r(\alpha)=s(\alpha)$, and we call the vertex
$s(\alpha)=r(\alpha)$ the \emph{base point} of the loop.  A loop is said to be
\emph{simple} if $s(\alpha_i) \neq s(\alpha_1)$ for all $1 < i \leq n$.  The
following is an important condition for graphs to satisfy.

$\text{ }$

\noindent \textbf{Condition~(K)}:  No vertex in $G$ is the
base point of exactly one simple loop; that is, every vertex
is either the base point of no loops or at least two simple loops.

$\text{ }$

The graph algebra $C^*(G)$ is unital if and only if $G$ has a
finite number of vertices, cf.~\cite[Proposition~1.4]{KPR}, and in
this case $1_{C^*(G)} = \sum_{v \in G^0} p_v$.  If $G$ has an
infinite number of vertices and we list them as
$G^0 = \{v_1, v_2, \ldots \}$ and define $p_n := \sum_{i=1}^n
p_{v_i}$, then $\{ p_n \}_{n=1}^\infty$ will be an approximate
unit for $C^*(G)$.

\begin{definition}  A \emph{trace} on a $C^*$-algebra $A$ is
a linear functional $\tau : A \rightarrow \C$ with the
property that $\tau(ab)=\tau(ba)$ for all $a,b \in A$.  We say
that $\tau$ is \emph{positive} if $\tau (a) \geq 0$ for all $a
\in A^+$.  If $\tau$ is a positive trace and $\| \tau \| = 1$ we call
$\tau$ a \emph{tracial state}.  The set of all tracial states is
denoted $T(A)$.
\end{definition}

\begin{definition}
If $G$ is a graph, then a \emph{graph trace} on $G$ is a
function $g : G^0 \rightarrow \R^+$ with the
following two properties: 
\begin{enumerate}
\item \label{g-t-1} For any $v \in G^0$ with $0 <
|s^{-1}(v)| < \infty$ we have $g(v) = \sum_{s(e) = v} g(r(e))$. 
\item \label{g-t-2} For any infinite emitter $v \in G^0$ and
any finite set of edges $e_1,
\ldots, e_n \in s^{-1}(v)$ we have $g(v) \geq
\sum_{i=1}^n g(r(e_i))$.
\end{enumerate}
\end{definition}

Because the value of $g$ at any vertex is non-negative, it
follows that whenever $v$ is an infinite emitter the infinite sum $\sum_{s(e)
= v} g(r(e))$ converges, and moreover
$\sum_{s(e) = v} g(r(e)) \leq g(v)$.

We define the \emph{norm} of $g$ to be the (possibly infinite) value $\|g
\| := \sum_{v \in G^0} g(v)$.  We shall call a graph trace
\emph{bounded} if $\| g \| < \infty$, and we shall use
$T(G)$ to denote the set of all graph traces on $G$ with
norm one.  Also note that if $v, w \in G^0$, then $v \geq w$ implies $g(v)
\geq g(w)$.

If $\tau : C^*(G) \to \C$ is a tracial state, then $\tau$ induces a graph
trace $g_\tau$ of norm one given by $g_\tau(v) := \tau(p_v)$.  If
$G$ satisfies Condition~(K), then the map $\tau \mapsto g_\tau$ is a bijection
(in fact, an affine homeomorphism) from $T(C^*(G))$ onto $T(G)$
\cite[\S3]{Tom6}.  There are examples which show that in general this
map is not injective.

\begin{definition}
We say that two projections $p,q \in A$ are equivalent, written $p \sim q$,
if there exists an element $v \in A$ with $p=vv^*$ and $q=v^*v$.
\end{definition}

In \cite{Cun4} Joachim Cuntz introduced a notion of comparison of (positive)
elements in a $C^*$-algebra for the purpose of constructing dimension
functions and traces on $C^*$-algebras.

\begin{definition}[Cuntz]
Let $A$ be a $C^*$-algebra, and let $a, b$ be positive elements
in $A$. We write $a \lesssim b$ if there exists a sequence
$\{x_k\}_{k=1}^\infty$ in $A$ with $x_k^* b x_k \to a$.
\end{definition}

If $p, q$ are projections in a $C^*$-algebra $A$, then $p \lesssim q$ if and
only if $p$ is equivalent to a subprojection of $q$; that is, there
exists a partial isometry $v \in A$ such that $p=vv^*$ and $v^*v \leq q$. 
Thus the above definition agrees with usual definition of comparison of two
projections.

If $e \in G^1$ then we see that $p_{r(e)} = s_es_e^*$ and $s_es_e^* \leq
p_{s(e)}$.  Therefore $p_{r(e)} \lesssim p_{s(e)}$.  More generally we see
that $v \geq w$ implies $p_w \lesssim p_v$.

\begin{definition}
If $G$ is a graph, a subset $H \subseteq G^0$ is said to be \emph{hereditary}
if for every $e \in G^1$ we have that $s(e) \in H$ implies $r(e) \in H$.  A
hereditary subset is said to be \emph{saturated} if whenever $v \in G^0$ with
$0 < | s^{-1}(v) | < \infty$ then $\{ r(e) : e \in G^1 \text{ and } s(e) = v
\} \subseteq H$ implies that $v \in H$.  If $H$ is a hereditary subset, then
the \emph{saturation} of $H$ is the smallest saturated hereditary subset
$\overline{H}$ of $G^0$ containing $H$. 
\end{definition}

If $H$ is a hereditary subset of $G^0$, then we can give an inductive
description of the saturation $\overline{H}$.  We define $H_0 := H$ and
having defined $H_n$ we set $$H_{n+1} := H_n \cup \{ v \in G^0 : 0 < |
s^{-1}(v) | < \infty \text{ and } s(e) = v \text{ implies } r(e) \in H_n
\}.$$ Then it is straightforward to show that $\overline{H} =
\bigcup_{n=0}^\infty H_n$.

\begin{definition}
Given a saturated hereditary subset $H \subseteq G^0$, we define 
$$ B_H := \{v \in G^0 : |s^{-1}(v)| = \infty \text{ and } 0 <
|s^{-1}(v) \cap r^{-1}(G^0 \setminus H)| < \infty \}.$$
Since $H$ is hereditary, we see that $B_H$ is disjoint from $H$.  If
$\{s_e, p_v\}$ is a generating Cuntz-Krieger $G$-family in $C^*(G)$, then for
$S \subseteq B_H$ we define
$$I_{(H,S)} := \text{ the ideal in $C^*(G)$ generated by $\{p_v : v
\in H \} \cup \{p_{v}^H : v \in S \}$},$$
where 
$$ p_{v}^H := p_{v} - \sum_{{s(e) = v} \atop {r(e) \notin H}} s_e s_e^*.$$
\end{definition}

\begin{definition}
If $H$ is a saturated hereditary subset of $G$ and $S \subseteq B_H$, then we
define a graph $G_{(H,S)}$ as follows:
\begin{align*}
G_{(H,S)}^0 &:= (G^0 \backslash H) \cup \{ v' : v \in B_H \backslash S \} \\
G_{(H,S)}^1 &:= \{e \in G^1 : r(e) \notin H \} \cup \{ e' : r(e) \in B_H
\backslash S \}
\end{align*}
and we extend $r$ and $s$ to $G_{(H,S)}^1$ by $r(e') = r(e)'$ and $s(e') =
s(e)$.  It follows from \cite[Corollary~3.5]{BHRS} that $C^*(G) / I_{(H,S)}
\cong C^*(G_{(H,S)})$. 
\end{definition}

\section{Stability of graph $C^*$-algebras} \label{stability-sec}

This section is devoted to proving Theorem~\ref{stable-thm}, which is a
generalization of \cite[Theorem~2.14]{Hje2}.

\begin{definition}  If $v$ is a vertex in a graph $G$ we define $L(v) := \{ w
\in G^0 : w \geq v \}$.  We say that $v$ is \emph{left infinite} if $L(v)$
contains infinitely many elements, and we say that $v$ is \emph{left finite}
if $L(v)$ contains finitely many elements.
\end{definition}

\begin{theorem} \label{stable-thm}
If $G$ is a graph, then the following are equivalent.
\begin{enumerate}
\item[(a)] $C^*(G)$ is stable
\item[(b)] $C^*(G)$ has no nonzero unital quotients and no tracial states
\item[(c)] Every vertex in $G$ that is on a loop is left infinite and
$T(G) = \emptyset$
\item[(d)] Every vertex in $G$ that is on a loop is left infinite and
$G$ has no nonzero bounded graph traces
\item[(e)] For every $v \in G^0$ and every finite set $F \subseteq G^0$ there
exists a finite set $W \subseteq G^0$ with $W \cap F = \emptyset$ and $p_v
\lesssim \sum_{w \in W} p_w$.
\item[(f)] For every finite set $V \subseteq G^0$ there exists a finite set
$W \subseteq G^0$ with $V \cap W = \emptyset$ and $\sum_{v \in V} p_v
\lesssim \sum_{w \in W} p_w$.
\end{enumerate}
\end{theorem}

\begin{corollary} \label{left-inf-stable}
If $G$ is a graph and every vertex of $G$ is left infinite, then $C^*(G)$ is
stable.  If $G$ has no sources and $C^*(G)$ is stable, then every vertex of
$G$ is left infinite.
\end{corollary}
\begin{proof}
Suppose every vertex of $G$ is left infinite.  If $v \in G^0$ and $F \subseteq
G^0$ is a finite set, then we may choose an element $w \in G^0$ such that $w
\notin F$ and $w \geq v$.  But then $p_v \lesssim p_w$ and by
Theorem~\ref{stable-thm}(e) $C^*(G)$ is stable.

If $G$ has no sources then for every $v \in G^0$ there exists a sequence of
edges $e_1e_2e_3 \ldots$ with $r(e_{i+1}) = s(e_i)$ and $r(e_1)=v$.  If the
elements of $\{s(e_i)\}_{i=1}^\infty$ are distinct, then $v$ is left
infinite.  If the elements of $\{s(e_i)\}_{i=1}^\infty$ are not distinct,
then there exists a loop that can reach $v$.  If $C^*(G)$ is stable, then by
Theorem~\ref{stable-thm}(c) all vertices on loops are left infinite.  Hence
$v$ is also left infinite.
\end{proof}

We cannot remove the condition of no sources in the converse of the above
corollary.  If $G$ is the graph 
$$
\xymatrix{
\bullet \ar[r] & \bullet \ar[r] & \bullet \ar[r] & \bullet \ar[r] & \ldots \\
}
$$
then no vertex of $G$ is left infinite, but $C^*(G) \cong \K$ is stable.

\begin{remark}
The equivalence of Conditions (a), (b), and (f) in Theorem~\ref{stable-thm}
was established for locally finite graphs in \cite[Theorem~2.14]{Hje2}.
We mention that Condition~(c) of Theorem~\ref{stable-thm} is
often easier to verify than Condition~(b).  This is because graph traces are
typically easier to deal with than tracial states, and it is often easy to
deduce whether $T(G)$ is empty simply by looking at $G$.  Furthermore, we
point out that the tracial states of $C^*(G)$ and the graph traces on $G$ of
norm one are not generally in one-to-one correspondence (see
\cite[\S3]{Tom6}).
\end{remark}

\begin{remark}  We see from Theorem~\ref{stable-thm} that a graph
$C^*$-algebra is stable if and only if it has no nonzero unital quotients and
no tracial states.  It is always the case that any stable $C^*$-algebra will
have no nonzero unital quotients and no tracial states, but in general the
converse does not hold.  (Interestingly, it is shown in
\cite[Proposition~5.1]{HR} that the converse will hold if certain full
hereditary subalgebras of the $C^*$-algebra satisfy a particular property.)
\end{remark}

\begin{lemma} \label{stab-equiv}
Let $A$ be a $C^*$-algebra with an increasing countable approximate unit $\{
p_n \}_{n=1}^\infty$ consisting of projections.  Then the following are
equivalent.
\begin{enumerate}
\item[(i)]  $A$ is stable.
\item[(ii)] For every projection $p \in A$ there exists a projection $q \in
A$ such that $p \sim q$ and $p \perp q$.
\item[(iii)] For all $n \in \N$ there exists $m > n$ such that $p_n \lesssim
p_m-p_n$ 
\end{enumerate}
\end{lemma}

\begin{proof}
The equivalence of \textrm{(i)} and \textrm{(ii)} is shown in
\cite[Theorem~3.3]{HR}.  The equivalence of \textrm{(ii)} and \textrm{(iii)}
is shown in \cite[Lemma~2.1]{Hje2}.
\end{proof}

\begin{lemma} \label{zero-sat}
If $G$ is a graph and $g:G^0 \to \R^+$ is a graph trace on $G$, then $$H := \{
v \in G^0 : g(v) = 0 \}$$ is a saturated hereditary subset.
\end{lemma}

\begin{proof}
If $e \in G^1$, then $g(s(e)) \geq g(r(e))$.  Thus $s(e) \in H$ implies that
$r(e) \in H$, and $H$ is hereditary.  If $v \in G^0$  is not a singular
vertex and $\{ r(e) : e \in G^1 \text{ and } s(e) = v \} \subseteq H$, then
$g(v) = \sum_{s(e)=v} g(r(e)) = 0$ so $v \in H$, and $H$ is saturated.
\end{proof}

\begin{lemma} \label{quotient-lift}
Let $G$ be a graph, let $H$ be a saturated hereditary subset of $G^0$, and
let $\pi : C^*(G) \to C^*(G) / I_{(H,\emptyset)}$ be the projection map.  If
$p$ is a projection in $C^*(G)$, $W \subseteq G^0 \backslash H$ is a finite
set, and $\pi(p) \lesssim \sum_{w \in W} \pi (p_w)$ in $C^*(G) /
I_{(H,\emptyset)}$, then there exists a finite set $X \subseteq H$ such
that $p \lesssim \sum_{w \in W} p_w + \sum_{x \in X} p_x$ in $C^*(G)$.
\end{lemma}

\begin{proof}
Write $H = \{ v_1, v_2, \ldots \}$.  If we let $p_n :=
\sum_{i=1}^n p_{v_i}$, then $I_{(H,\emptyset)}$ is generated by $\mathcal{P}
= \{p_n\}_{n=1}^\infty$ and \cite[Lemma~2.6]{Hje2} implies that $p \lesssim
\sum_{w \in W} p_w + p_n$ for some $n$.
\end{proof}

\noindent \emph{Proof of Theorem~\ref{stable-thm}.}  (a) $\Longrightarrow$ (b)
:  It is shown in \cite[Proposition~5.1]{HR} that stable $C^*$-algebras have
no nonzero unital quotients and admit no nonzero traces.
\\

\noindent (b) $\Longrightarrow$ (c) :  We shall first show that every vertex
on a loop is left infinite.   Let $\alpha$ be a loop in $G$ that is
based at $v$.  Then $H := G^0 \backslash L(s(\alpha))$ is a saturated
hereditary subset.  By hypothesis $C^*(G) / I_{(H,B_{H})} \cong G_{(H,B_{H})}$
is nonunital and hence $G_{(H,B_{H})}^0 = G^0 \backslash H = L(s(\alpha))$ is
infinite.  Thus $s(\alpha)=v$ is left infinite. 

Now we shall show that $T(G)$ is empty by supposing that there exists $g
\in T(G)$ and arriving at a contradiction.  Let us begin by showing that if
$v$ is a vertex on a loop, then $g(v) = 0$.   From the previous paragraph
every vertex on a loop is left infinite.  Since $$\|g\| = \sum_{w \in G^0}
g(w) \geq \sum_{w \in L(v)}g(w)$$ and since $w \geq v$ implies $g(w) \geq
g(v)$ the only way that this infinite sum can be finite is if $g(v)=0$.  Thus
$g$ vanishes on every vertex that is on a loop.

If we now let $H := \{ v \in G^0 : g(v) = 0 \}$ then it follows from
Lemma~\ref{zero-sat} that $H$ is a saturated hereditary subset.  We define a
graph trace $\tilde{g}$ on $G_{(H,\emptyset)}$ by $$\tilde{g} (w) :=
\begin{cases} g(w) & \text{ if $w \in (G^0 \backslash H) \backslash
B_{H}$} \\ & \\ \displaystyle \sum_{ { s(e)=w } \atop {r(e) \notin H} }
g(r(e)) & \text{ if $w \in B_{H}$} \\ g(v) - \displaystyle \sum_{ { s(e)=v}
\atop {r(e) \notin H} } g(r(e)) & \text{ if $w = v'$ for some
$v \in B_{H}$.}
\end{cases}$$  It is straightforward to verify that $\tilde{g}$ is a graph
trace on $G_{(H,\emptyset)}$ and that $\| \tilde{g} \| = \| g \| = 1$.  Now
it follows from the previous paragraph that there are no loops in $G$ with
vertices in $G^0 \backslash H$.  Hence $G_{(H,\emptyset)}$ is a graph with
no loops.  Therefore \cite[\S3.3]{Tom6} implies that there exists a
tracial state $\tau$ on $C^*(G_{(H,\emptyset)})$.  Since $C^*(G) /
I_{(H,\emptyset)} \cong C^*(G_{(H,\emptyset)})$ it follows that $\tau$ lifts
to a tracial state on $C^*(G)$.  But this contradicts the fact that $C^*(G)$
has no tracial states.
\\

\noindent (c) $\Longrightarrow$ (d) : If $g$ was a nonzero bounded graph
trace on $G$, then we could normalize to get an element $\frac{1}{\|g\|} \cdot
g \in T(G)$.
\\

\noindent (d) $\Longrightarrow$ (e) :  Choose a vertex $v \in G^0$.  Define
$H := \{ w \in G^0 : \text{ $w$ is left infinite} \}$.  Then $H$ is a
hereditary subset, and we let $\overline{H}$ denote the saturation of $H$. 
Consider the following two cases:

\noindent \textsc{Case I:} $v \in \overline{H}$.  Define $H_0 := H$ and for
each $n \in \N$ set $$H_{n+1} := H_n \cup \{ w \in G^0 : 0 < | s^{-1}(w) | <
\infty \text{ and } s(e) = w \text{ implies } r(e) \in H_n \}.$$ Then we see
that $\overline{H} = \bigcup_{n=0}^\infty H_n$. We shall prove that the claim
holds whenever $v \in \overline{H}$ by induction on $k := \min \{ n \in \N : v
\in H_n \}$.  In the base case we have $k=0$ and thus $v \in H$.  Since every
vertex in $H$ is left infinite for every finite set $F \subseteq G^0$ there
exists $w \in G^0$ such that $w \notin F$ and $w \geq v$.  But then $p_v
\lesssim p_w$ and the claim holds.  Now assume that the claim holds whenever
$v$ is in $\overline{H}$ with $\min \{ n \in \N : v \in H_n \}$ strictly less
than a fixed $k$.  Suppose that $v \in H_k$.  Then $s^{-1}(v)$ consists of a
finite and nonzero number of edges $\{ e_1, \ldots, e_n \}$ with $r(e_i)
\in H_{k-1}$ for all $i$.  By the induction hypothesis there exists a finite
set $W_1 \subseteq G^0$ such that $W_1 \cap F = \emptyset$ and $p_{r(e_1)}
\lesssim \sum_{w \in W_1} p_w$.  Similarly for each $1 < i \leq n$ there
exists a finite set $W_i \subseteq G^0$ which is disjoint from $F \cup W_1
\cup \ldots \cup W_{i-1}$ and with $p_{r(e_i)} \lesssim \sum_{w \in W_i}
p_w$.  Now if we let $x = s_{e_1} + \ldots + s_{e_n}$ then we see that $x^*
(\sum_{i=1}^n s_{e_i}s_{e_i}^*) x = \sum_{i=1}^n p_{r(e_i)}$, and thus
$\sum_{i=1}^n s_{e_i}s_{e_i}^* \lesssim \sum_{i=1}^n p_{r(e_i)}$.  Therefore
if we let $W := W_1 \cup \ldots \cup W_n$ we see that
$W \cap F = \emptyset$ and $p_v = \sum_{i=1}^ns_{e_i}s_{e_i}^* \lesssim 
\sum_{i=1}^n p_{r(e_i)} \lesssim \sum_{w \in W_1} p_w + \ldots \sum_{w \in
W_n} p_w = \sum_{w \in W} p_w$.

\noindent \textsc{Case II:} $v \notin \overline{H}$.  Since every vertex on a
loop is left infinite, it follows that no vertices of $G^0 \backslash
\overline{H}$ are on loops.  Thus $G_{(\overline{H}, \emptyset)}$ contains no
loops and \cite[Corollary~2.13]{DT1} implies that $C^*(G_{(\overline{H},
\emptyset)})$ is an AF-algebra.  Furthermore, there are no tracial
states on $C^*(G_{(\overline{H}, \emptyset)}) \cong C^*(G) /
I_{(\overline{H},\emptyset)}$ since any tracial state would lift to a tracial
state on $C^*(G)$ and thus induce a graph trace of norm one on $G$.  Since
$C^*(G_{(\overline{H}, \emptyset)})$ is an AF-algebra with no tracial
states it follows from \cite[Theorem~4.10]{Bla4} that it is stable.  

If we list the vertices of $G^0 \backslash \overline{H}$ as $\{ w_1, w_2,
\ldots \}$ with $w_1 = v$, then the elements $p_n := \sum_{i=1}^n \pi
(p_{w_i})$ form an increasing approximate unit for $C^*(G_{(\overline{H},
\emptyset)})$ consisting of projections.  If $F \subseteq G^0$ is a finite
set, let $n = \max \{ i \in \N : w_i \in F \}$.  Since
$C^*(G_{(H,\emptyset)})$ is stable, Lemma~\ref{stab-equiv}(c) implies that
there exists $m > n$ such that $p_n \lesssim p_m-p_n$.  But if we let $W_0 :=
\{ w_{n+1}, \ldots, w_m\}$ then $W_0 \cap F = \emptyset$ and $\pi(p_v)
\lesssim p_n \lesssim p_m-p_n = \sum_{w \in W_0} \pi(p_w)$.  It then follows
from Lemma~\ref{quotient-lift} that there exists a finite set $X \subseteq
\overline{H}$ for which $p_v \lesssim \sum_{w \in W} p_w + \sum_{x \in X} p_x$
in $C^*(G)$.  Now since $X \subseteq \overline{H}$ we see from Case~I above
that if $X = \{x_1 , \ldots, x_n \}$ then for each $i$ we may choose $W_i$
such that $W_i$ is disjoint from $F \cup W_0 \cup \ldots W_{i-1}$ and $p_{x_i}
\lesssim \sum_{w \in W_i} p_w$.  If we let $W := W_0 \cup \ldots \cup W_n$,
then $W \cap F = \emptyset$ and $p_v \lesssim \sum_{w \in W_0}p_w + \sum_{x
\in X}p_x \lesssim \sum_{w \in W_0}p_w + \ldots + \sum_{w \in W_n}p_w =
\sum_{w \in W} p_w$.
\\

\noindent (e) $\Longrightarrow$ (f) : List the elements of $V$ as $V = \{
v_1, \ldots v_n \}$.  Choose $W_1$ such that $W_1 \cap V = \emptyset$ and
$p_{v_1} \lesssim \sum_{w \in W_1} p_w$.  Having chosen $W_k$ we may choose
$W_{k+1}$ so that $W_{k+1}$ is disjoint from $V \cup W_1 \cup \ldots \cup W_k$
and $p_{v_k} \lesssim \sum_{w \in W_k} p_w$.  We continue in
this fashion until we produce $n$ sets $W_1, \ldots, W_n$ with these
properties.  If we let $W := W_1 \cup \ldots \cup W_n$, then $V \cap W =
\emptyset$ and  $\sum_{v \in V} p_v \lesssim \sum_{w \in W_1} p_w + \ldots +
\sum_{w \in W_n} p_w = \sum_{w \in W} p_w$.
\\

\noindent (f) $\Longrightarrow$ (a) :  List the vertices of $G$ as $G^0 := \{
v_1, v_2, \ldots \}$.  For each $n \in \N$ we define $p_n := \sum_{i=1}^n
p_{v_i}$.  Then $\{p_n \}_{n=1}^\infty$ is an increasing approximate unit
consisting of projections, and by Lemma~\ref{stab-equiv} it suffices to prove
that for all $n \in N$ there exists $m > n$ such that $p_n \lesssim
p_m-p_n$.

Let $n \in N$, and define $V := \{ v_1, \ldots, v_n \}$.  By hypothesis there
exists a finite set $W \subseteq G^0$ such that $V \cap W = \emptyset$ and
$\sum_{v \in V} p_v \lesssim \sum_{w \in W} p_w$.  Let $m:= \max \{ k \in \N :
v_k \in W \}$.  Since $V \cap W = \emptyset$ we see that $\sum_{w \in W} p_v
\leq p_m - p_n$.  Thus $p_n = \sum_{v \in V} p_v \lesssim \sum_{v \in W} p_v
\leq p_m-p_n$.
\hfill $\qed$

\section{The stabilization of a graph $C^*$-algebra} \label{stabilization-sec}

\begin{definition}
If $G$ is a graph and $v \in G^0$ is a vertex, then by \emph{adding a head to
$v$} we mean attaching a graph of the form 
$$
\xymatrix{
\cdots \ar[r]^{e_4} & v_3 \ar[r]^{e_3} & v_2 \ar[r]^{e_2} & v_1 \ar[r]^{e_1} &
v\\ }
$$
Thus we create a new graph $F$ from $G$ by defining $F^0 := G^0 \cup \{ v_1,
v_2, \ldots \}$, $F^1 := G^1 \cup \{e_1, e_2, \ldots \}$, and extend $r$ and
$s$ to $F^1$ by $r(e_i) = v_{i-1}$ and $s(e_i) = v_i$.

The terminology ``adding a head" is meant to complement the
terminology for the analogous concept of ``adding a tail" introduced in
\cite[(1.2)]{BPRS}.
\end{definition}

\begin{theorem} \label{stabilization-gr-alg-thm}
If $G$ is a graph, let $\tilde{G}$ be the graph obtained by adding a head
to each vertex of $G$.  Then $C^*(\tilde{G})$ is the stabilization of
$C^*(G)$; that is, $$C^*(\tilde{G}) \cong C^*(G) \otimes \K.$$
\end{theorem}

\begin{proof}
Following the proof of \cite[Lemma~1.2]{BPRS} one can show that $C^*(G)$ is
naturally isomorphic to a full corner of $C^*(\tilde{G})$.  Consequently
$C^*(G)$ is Morita equivalent to $C^*(\tilde{G})$, and since
Corollary~\ref{left-inf-stable} implies that $C^*(\tilde{G})$ is stable we
have $C^*(G) \otimes \K \cong C^*(\tilde{G}) \otimes \K \cong C^*(\tilde{G})$.
\end{proof}

\begin{corollary}
The class of graph $C^*$-algebras is closed under stabilization.
\end{corollary}

\begin{example}
If $G$ is the graph 
$$
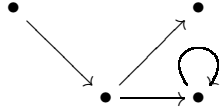
\xymatrix{
\bullet \ar[rd] & & \bullet \\ & \bullet \ar[ru] \ar[r] & \bullet \ar@(ul,ur)
}
$$
then $\tilde{G}$ is the graph
$$
\xymatrix{
\cdots \ar[r] & \bullet \ar[r] & \bullet \ar[r] & \bullet \ar[rd] & & \bullet 
& \bullet \ar[l] & \bullet \ar[l] & \cdots \ar[l] \\ 
\cdots \ar[r] & \bullet \ar[r] & \bullet \ar[r] & \bullet \ar[r] & \bullet
\ar[ru] \ar[r] & \bullet \ar@(ul,ur) & \bullet \ar[l] & \bullet \ar[l] &
\cdots \ar[l] }
$$
and $C^*(\tilde{G}) \cong C^*(G) \otimes \K$.
\end{example}

\begin{example}
If $G$ is the following graph with one vertex and infinitely many edges, then
$C^*(G) \cong \mathcal{O}_\infty$
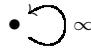
$$
\xymatrix{
\bullet \ar@(dr,ur)_\infty}
$$
and $\tilde{G}$ is the graph
$$
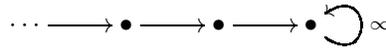
\xymatrix{
\cdots \ar[r] & \bullet \ar[r] & \bullet \ar[r] & \bullet \ar@(dr,ur)_\infty }
$$
so that $C^*(\tilde{G}) \cong \mathcal{O}_\infty \otimes \K$.
\end{example}

\begin{remark}
To obtain the stabilization it is often unnecessary to add a head to every
vertex in $G$.  It suffices to add enough heads to make all vertices
left infinite.  For example, one could choose to add heads only at the left
finite vertices of
$G$.
\end{remark}

\end{document}